\documentclass{amsart}
\usepackage{amscd,amssymb}
\theoremstyle{plain}
\newtheorem{thm}{Theorem}[section]
\newtheorem{lem}[thm]{Lemma}

\newtheorem{prop}[thm]{Proposition}

\theoremstyle{definition}
\newtheorem{Def}[thm]{Definition}

\theoremstyle{remark}

\newtheorem{Remark}[thm]{Remark}

\numberwithin{equation}{section}

\usepackage{pst-all}

\begin{document}
\setcounter{section}{-1}

\title{Non-nesting actions of Polish groups on real trees}
\author{{Vincent Guirardel, Aleksander Ivanov}}
\thanks{The second author is supported by KBN grant 2 P03A 007 19}
\thanks{{\em E-mail addresses}:  vincent.guirardel@math.univ-toulouse.fr {\em and}  ivanov@math.uni.wroc.pl}
\thanks{{\em Fax number}: 48-71-3757429 (A.Ivanov)}

\maketitle 
\centerline{Institut de Math\'ematiques de Toulouse, Universit\'e Paul Sabatier Toulouse 3}
\centerline{31062 Toulouse cedex 9, France} 
\centerline{Institute of Mathematics, Wroc{\l}aw University, pl.Grunwaldzki 2/4, }
\centerline{50-384, Wroc{\l}aw, Poland}
\bigskip
{\bf Abstract}
\begin{quote}
We prove that if a Polish group $G$ with a comeagre conjugacy 
class has a non-nesting action on an $\mathbb{R}$-tree, then 
every element of $G$ fixes a point.
\parskip0pt


\end{quote}

\section{Introduction}

Non-nesting actions by homeomorphisms on $\mathbb{R}$-trees 
frequently arise in geometric group theory.
For instance, they occur in Bowditch's study of cut points of the boundary
at infinity of a hyperbolic group \cite{bow}, or in the Drutu-Sapir study of tree-graded spaces \cite{DS},
and their relations with isometric actions were studied in \cite{levitt}. 
Non-nesting property is a topological substitute for an isometric action. 
It asks that no interval of the $\mathbb{R}$-tree is sent properly into itself by an element of the group.
\parskip0pt 

In this paper, we are concerned with a Polish group $G$ having a comeagre conjugacy class.
The group $S_{\infty}$ of all permutations of $\mathbb{N}$ and more generally 
the automorphism group of any $\omega$-stable $\omega$-categorical 
structure (see \cite{hhls}) provide typical model-theoretic examples.  
Among other examples, we mention the automorphism group of the random graph 
and the groups $Aut(\mathbb{Q}, <)$, $Homeo (2^{\mathbb{N}})$ and $Homeo_{+}(\mathbb{R})$.  
The latter ones appear in \cite{kecros} and \cite{solros} as important 
cases of extreme amenability and automatic continuity of homomorphisms.   
The property of having a comeagre conjugacy class 
plays an essential role in these respects.  

The following theorem is the main result of the paper: 
\begin{quote}
Consider a group $G$ with a non-nesting action on an $\mathbb{R}$-tree $T$.
If $G$ is a Polish group with a comeagre conjugacy class,
then every element of $G$ fixes a point in $T$.
\end{quote}
This theorem generalizes the main result of the paper of 
H.D.Macpherson and S.Thomas \cite{mactho} (where the authors 
study actions of Polish groups on simplicial trees) and 
extends Section 8 of the paper of Ch.Rosendal \cite{rosendal} 
(concerning isometric actions on $\Lambda$-trees). 
It is worth noting that some related problems 
have been studied before 
(see \cite{bow}, \cite{chi}, \cite{dun} and \cite{levitt}).
Our motivation is partially based on these investigations.

\section{Non-nesting actions on $\mathbb{R}$-trees}

\begin{Def}
An $\mathbb{R}$-tree is a metric space $T$ such that for any $x\neq y\in T$,
there is a unique topologically embedded arc joining $x$ to $y$, and this arc 
is isometric to some interval of $\mathbb{R}$.
\end{Def} 

Equivalently, as a topological space, $T$ is a metrizable, uniquely arc-connected, 
locally arc-connected topological space \cite{MaOv}.
We define $[x,y]$ as the arc joining $x$ to $y$ if $x\neq y$, and $[x,y]=\{x\}$ if $x=y$.
We say that $[x,y]$ is a \emph{segment}. 

A subset $S\subseteq T$ is \emph{convex} if 
$(\forall x,y \in S) [x,y] \subseteq S$. 
A convex subset is also called a \emph{subtree}. 
Given $x,y,z \in T$, 
there is a unique element  $c \in [x,y]\cap [y,z] \cap [z,x]$,
called the {\em median} of $x,y,z$. 
When $c\notin\{x,y,z\}$, the subtree $[x,y]\cup [x,z]\cup [y,z]$ is called a \emph{tripod}.
A \emph{line} is a convex subset containing no tripod and maximal for inclusion.

Given two disjoint closed subtrees $A,B\subseteq T$, there exists 
a unique pair of points $a\in A,b\in B$ such that for all $x\in A,y\in B$, 
$[x,y]\supseteq [a,b]$. 
The segment $[a,b]$ is called the \emph{bridge} between $A$ and $B$.
If $x\notin A$, the \emph{projection} of $x$ on $A$ is the point $a\in A$ 
such that $[x,a]$ is the bridge between $\{x\}$ and $A$.
\\

The betweenness relation $B$ of $T$ is the ternary relation
$B(x;y,z)$ defined by $x\in (y,z)$.
A \emph{weak homeomorphism} of the $\mathbb{R}$-tree $T$
is a bijection $g:T\to T$ which preserves the betweenness relation. 
Any homeomorphism of $T$ is clearly a weak homeomorphism.
All actions on $T$ are via weak homeomorphisms.

\begin{Remark}
If $g:T\rightarrow T$ is a weak homeomorphism, then its restriction 
to each segment, to each line, and to each finite union of segments
is a homeomorphism onto its image (for the topology induced by the metric). 
This is because the metric topology agrees with the
topology induced by the order on a line or a segment.
Conversely, any bijection $g:T\to T$ which maps each segment homeomorphically onto 
its image is a weak homeomorphism as it maps $[x,y]$ to the unique embedded arc joining $g(x)$ to $g(y)$.
\end{Remark}

\begin{Remark}\label{rem_closed}
If $S\subseteq T$ is a subtree, then $S$ is closed (for the topology induced by the metric)
if and only if $S\cap I$ is closed in $I$ for every segment $I$. 
In particular, a weak homeomorphism preserves the set of closed subtrees.  
\end{Remark}

\begin{Def}
An action of $G$ on $T$ by weak homeomorphisms is \emph{non-nesting} if there 
is no segment $I\subseteq T$, and no $g\in G$ such that $g(I)\subsetneqq I$. 
\end{Def}

From now on, we assume that $G$ has a non-nesting action on an $\mathbb{R}$-tree $T$.
We say that $g\in G$ is \emph{elliptic} if it has a fixed point, and \emph{loxodromic} otherwise.

\begin{lem}[{\cite[Theorem 3]{levitt}}]\label{lem_nnest}
Let $G$ be a group with a non-nesting action on an $\mathbb{R}$-tree $T$.

  \begin{itemize}
  \item If $g$ is elliptic, its set of fix points $T^g$ is a closed convex subset.
  \item\label{nnest_it2} If $g$ is loxodromic, there exists a unique line $L_g$ preserved by $g$; moreover,
$g$ acts on $L_{g}$ by  an order preserving transformation, which is a translation 
  up to topological conjugacy. 
  \end{itemize} 
\end{lem}

In \cite{levitt}, $g$ is assumed to be a homeomorphism, but the argument still applies,
except to prove that $T^g$ is closed.
This fact  follows from Remark \ref{rem_closed}.

When $g$ is loxodromic, $L_g$ is called the \emph{axis} of $g$.
The action of $g$ on $L_g$ defines a natural ordering on $L_g$ such that for all $x\in L_g$, $x<g(x)$.

The proof of the following lemma is standard
(by arguments from \cite{tits}, Section 3.1)
and can be found in \cite{ivanov}.

\begin{lem} \label{loxo}
If $g$ is loxodromic, then for any $p \in T$, $[p,g(p)]$ meets $L_{g}$ 
and $[p,g(p)] \cap  L_{g} = [q,g(q)]$ for some $q \in L_{g}$. 
\end{lem}

\begin{prop} \label{3.1}
Let $G$ be a group with a non-nesting action on an $\mathbb{R}$-tree $T$.
Then

\begin{enumerate}
\item \label{it1} If $g$ is elliptic and $x\notin T^g$, then $[x,g(x)]\cap T^g=\{a\}$ where 
$a$ is the projection of $x$ on $T^g$. 

\item\label{it2}
If $g,h\in G$ are elliptic and $T^{g}\cap T^h=\emptyset$,
then $gh$ is loxodromic, its axis contains the bridge between $T^g$ and $T^h$,
and $T^g\cap L_{gh}$ (resp.\ $T^h\cap L_{gh}$)
contains exactly one point.  
In particular, if $g,h$ and $gh$ are elliptic, then $T^g\cap T^h\cap T^{gh}\neq\emptyset$. 

\item \label{it4}
Let $h,h'\in G$ be loxodromic elements, and $a\in L_h$ be such that
for some $a'\in T$, $[a',(h')^2 (a')]\subseteq [a,h (a)]$. 
Then $h$ and $h'$ are not conjugate. 
\end{enumerate}
\end{prop}

These facts are classical for isometries of an $\mathbb{R}$-tree.
Assertion (\ref{it4}) is some substitute for the fact that the translation length of 
an isometry is a conjugacy invariant.

{\em Proof.} 
To prove Assertion (\ref{it1}),
consider $x\notin T^g$, and $I=[x,a]$ the bridge between $\{x\}$ and $T^g$.
If $g(I)\cap I=\{a\}$, we are done. 
Assume otherwise that $g(I)\cap I=[a,b]$ for some $b\neq a$.
Since $g(b)\neq b$, either $g.[a,b]\subsetneqq [a,b]$ or
$g.[a,b]\supsetneqq [a,b]$, in contradiction with the non-nesting assumption.

To see (\ref{it2}), consider $I=[a,b]$ the bridge between 
$T^{g}$ and $T^{h}$ with $a\in T^g,b\in T^h$, and let $J=h^{-1} (I)\cup I$.
By Assertion (\ref{it1}),  $I\cap h^{-1}(I)=\{b\}$ (resp.\ $I\cap g(I)=\{a\})$, $I\cap h(I)=\{b\})$,)
so $h^{-1}(a),b,a$ (resp. $b,a,g(b)$, $\  a,b,h(a)$ hence $a=g(a),g(b),gh(a)$) are aligned in this order.
In particular $h^{-1}(a),b,a,g(b),gh(a)$ are aligned in this order so  
$h^{-1} (I), I , g(I) , gh(I)$ are four consecutive non-degenerate subsegments 
of the segment $[h^{-1}(a),gh(a)]$.
This implies that $gh(J)\cap J=\{a\}$.
If $gh$ was elliptic, $J=[h^{-1}(a),gh(h^{-1}(a)]$ would contain a  point fixed by $gh$,
and this fix point would have to lie in $gh(J)\cap J$, but this is impossible since $gh(a)\neq a$.
We claim that $J\subseteq L_{gh}$. 
Otherwise, the segment $J_0=J\cap L_{gh}$ is a  proper subsegment of $J$,
and $gh(J_0)\cap J_0=\emptyset$, contradicting Lemma \ref{loxo}. 
Since $J\cap T^h=\{b\}$ and since $T^h$ is convex, $L_{gh}\cap T^h=\{b\}$.
Similarly, $(I\cup g(I))\cap T^g=\{a\}$ implies that $L_{gh}\cap T^g=\{a\}$.

Statement (\ref{it4}) is easy: let $I=[a,h (a)]\subseteq L_h$, and let $I'=[a',(h')^2 (a')]\subseteq I$. 
By Lemma \ref{loxo}, changing $I'$ to some subsegment, we may assume that 
$I'\subseteq L_{h'}$ so that $I'$ is a fundamental domain for the action of $(h')^2$ 
on $L_{h'}$ by Lemma \ref{lem_nnest}.
If $h' = h^g$, $g^{-1}(L_{h})=L_{h'}$ and $g^{-1}(I)$ is a fundamental domain 
for the action of $h'$ on $L_{h'}$.
Replacing $g$ by some $g (h')^i$ ($i\in \mathbb{Z}$), 
if necessary we obtain $g^{-1}(I) \subsetneqq I'\subseteq I$, 
a contradiction with the non-nesting assumption. 
$\square$

\section{Polish groups with comeagre conjugacy classes}

A Polish group is a topological group whose topology is
{\em Polish} (a Polish space is a separable completely
metrizable topological space).
A subset of a Polish space is {\em comeagre} if it contains
an intersection of a countable family of dense open sets.

H.D.Macpherson and S.Thomas have proved in \cite{mactho}
that if a Polish group has a comeagre conjugacy class then
every element of the group fixes a point under any action on
a $\mathbb{Z}$-tree without inversions. 
Ch.Rosendal has generalized this theorem to the case when 
the group acts on an $\Lambda$-tree by isometries 
(see Section 8 in \cite{rosendal}). 
In this section we consider the case of non-nesting actions.

\begin{thm} \label{4.2}
Consider a group $G$ with a non-nesting action on an $\mathbb{R}$-tree $T$.
If $G$ is a Polish group with a comeagre conjugacy class, 
 then every element of $G$ is elliptic.  
\end{thm}

\begin{Remark}
  We don't assume any relation between the action of $G$ and its topology as a Polish group:
the action of $g$ is not assumed to depend continuously on $g$.
\end{Remark}

\begin{Remark} 
  Using Proposition \ref{3.1}(\ref{it2}), one can extend the proof of 
Serre's Lemma \cite[Prop 6.5.2]{serre},
and show that every finitely generated subgroup of $G$ fixes a point in $T$.
It follows that $G$ fixes a point or an end of $T$.
\end{Remark}

We start with the following lemma. 

\begin{lem} \label{MaLe}
Under the circumstances of Theorem \ref{4.2}, assume that 
$h_1,h_2 \in G$ are conjugate and loxodromic, and that
$g=h_2 h_1$ is conjugate to $h^6_1$ or $h^{-6}_1$. 
Then $L_{h_1} \cap L_{h_2} =\emptyset$. \parskip0pt 

Moreover, denoting by $[a,b]$ the bridge between $L_{h_1}$ and $L_{h_2}$ 
with $a\in L_{h_1}$, $b\in L_{h_2}$ then 
$$
[h_1^{-1}(a),a]\cup [a,b]\cup [b,h_2 (b)] \subseteq  L_{g} 
$$
and  $h_1^{-1}(a) < a <b <h_2 (b)$ 
for the ordering of $L_g$ defined after Lemma \ref{lem_nnest}.
\end{lem} 

{\em Proof.} 
Assuming the contrary, consider $t\in L_{h_1} \cap L_{h_2}$ and 
$p=h^{-1}_1 (t)$. 
Since $[p,g(p)] \subseteq [h^{-1}_1 (t),t]\cup [t,h_2 (t)]$, 
may find  $q\in L_g$ such that  
$[q,g(q)]\subseteq [h^{-1}_1 (t),t]\cup [t,h_2 (t)]$.  

Consider $g_0$ such that $g_0^6=g$, and $g_0$ conjugate to $h_1$ or $h^{-1}_1$. 
Let $I=[q,g^2_0(q)]$. 
Since $L_{g_0}=L_g$, $I\subseteq L_{g_0}$ and
$I\cup g^2_0 (I) \cup g^4_0 (I)=[q,g_0^6(q)]\subseteq [h^{-1}_1 (t),t]\cup [t,h_2 (t)]$. 
Either $I$ or $g^4_0 (I)$ is contained in $[h^{-1}_1 (t),t]$ or in 
$[t,h_2 (t)]$, say $I\subseteq [h^{-1}_1 (t),t]$ for instance. 
Since $t\in L_{h_1}$, this contradicts Proposition \ref{3.1}(\ref{it4}). 

To see the final statement note that $L_g$ intersects 
 $[h^{-1}_1 h^{-1}_2 (a), a]$ and $[b, h_2 h_1(b)]$,
hence contains the bridge between these segments, i.e.\ $[a,b]$. 
It follows that $L_g$ contains $[h_1^{-1}h_2^{-1}(a),a]\supseteq[h_1^{-1}(a),a]$
and $[b,h_2h_1 (b)]\supseteq[b,h_2 (b)]$.
The lemma follows. 
$\square$
\bigskip

{\em Proof of Theorem \ref{4.2} }
Let $X$ be a conjugacy class of $G$ which is comeagre in $G$. 
Then $X\cap X^{-1}\neq \emptyset$,
but since $X$ is a conjugacy class $X=X^{-1}$. 
Note that 
\begin{quote} 
(*) For every sequence $g_{1},...,g_{m} \in G$ there exist
$h_{0},h_{1},...,h_{m} \in X$ such that for every $1 \le i \le m$,
$g_{i} = h_{0}h_{i}$. 
\end{quote} 

Indeed, let $g_{1},...,g_{m} \in G$.
Since $X$ and $g_{i}X^{-1}$ are comeagre in  $G$, all $g_{i}X^{-1}$ and $X$ 
have a common element $h_{0}\in X$.
Now there are $h_{1},...,h_{m} \in X$ such that for
any $1 \le i \le m$, $g_{i} = h_{0}h_{i}$.

First assume that $X$ consists of loxodromic elements, and argue towards a contradiction. 
Take $h\in X$ and consider $g=h^6$. 
By (*) above find $h_0 ,h_1 ,h_2 \in X$ such that 
$g=h_0 h_1$ and $g^{-1}=h_0 h_2$. 

Applying Lemma \ref{MaLe} to $h_0$, $h_1$ and to $h_0$, $h_2$, 
we get that $L_{h_0} \cap L_{h_1}=\emptyset$ and 
$L_{h_0} \cap L_{h_2}=\emptyset$. 
Let $b\in L_{h_0}$ and $a\in L_{h_1}$ define the 
bridge between $L_{h_0}$ and $L_{h_1}$,  
and let $b'\in L_{h_0}$ and $a'\in L_{h_2}$ define the 
bridge between $L_{h_0}$ and $L_{h_2}$. 
Since $L_g = L_{g^{-1}}$, by Lemma \ref{MaLe} we see 
that the segments $[a,b]\cup [b,h_0 (b)]$ and 
$[a',b']\cup [b',h_0 (b')]$ belong to $L_g$. 
Since $L_g$ does not contain a tripod, $b=b'$. 
Then $b<h_0 (b)$ both with respect to the order 
defined by $g$ and by $g^{-1}$. 
This is a contradiction,  so $X$ 
 consists of elliptic elements.

Assume that some $g \in G$ is loxodromic, and argue towards a contradiction.
Write $g = h'\cdot h$ for some $h,h' \in X$.
Then $T^h\cap T^{h'}=\emptyset$ and
denote by $I$ the bridge between  $T^{h}$ and $T^{h'}$. 
By Lemma \ref{3.1}(\ref{it2}) $I\subseteq L_g$.\parskip0pt

By (*) there exist  $h_{0},h_{1},h_{2},h_{3} \in X$ such that
$h = h_{0}h_{1}$, $h' = h_{0}h_{2}$, and $g = h_{0}h_{3}$.
By Lemma \ref{3.1}(\ref{it2}) there are
$a_{1} \in T^{h_{0}} \cap T^{h}$
and $b_{1} \in T^{h_{0}} \cap T^{h'}$. 
Then $I\subseteq [a_1,b_1]\subseteq T^{h_0}$. 
On the other hand, by Lemma \ref{3.1}(\ref{it2}) applied to $h_{0}$ and 
$h_{3}$, the intersection $T^{h_{0}} \cap L_{g}$ is a singleton.
Since $I$ is contained in this intersection, this is a contradiction. 
$\square$


\end{document}